\documentclass[11pt]{article}
\usepackage[T1]{fontenc}
\usepackage[utf8]{inputenc}
\usepackage{lmodern}
\usepackage[a4paper,margin=27mm]{geometry}
\usepackage{amsmath,amssymb,amsthm}
\usepackage{microtype}
\usepackage{xurl}
\usepackage[hidelinks,hypertexnames=false]{hyperref}
\title{Recovering Dirichlet coefficients from even samples}
\author{Payam Danesh}
\date{}
\begin{document}
\maketitle

\begin{abstract}
In this work, we determine the exact case error for recovering a
Dirichlet coefficient from finitely many positive even samples with
bounded additive errors. For sequences with a bounded tail and an
unrestricted initial block, a polynomial estimator is minimax among all
recovery maps. Its tail and noise constants admit explicit product
formulas. We prove a sharp tail asymptotic on the index scale
\(K^{3/2}\) and give a terminating procedure for selecting the sample
count from a computable error bound.

\end{abstract}


\hypertarget{introduction}{%
\section{Introduction}\label{introduction}}

First, for \(\mathbb{N} = \{1,2,\ldots\},\) we let \(q,K\in \mathbb{N}\)
and fix \(n \in \{ 1,\ldots,K\}.\) Here \(K\) is the sample count,
\(2q\) is the first sampling point and \(n\) is the target index. For a
real constant \(A > 0,\) we define

\begin{equation}
\mathcal{A}_{K}(A) = \{ a:\mathbb{N} \rightarrow \mathbb{C}:\ \left| a(k) \right| \leq A\text{ for every }k > K\}.
\tag{1.1}\label{eq:1.1}
\end{equation}

The first \(K\) coefficients are arbitrary complex numbers. The absence
of a common bound on this block is essential to the exact recovery
result.

For \(a \in \mathcal{A}_{K}(A),\) we define the ordinary Dirichlet
series by:

\begin{equation}
F_{a}(s) = \sum_{k = 1}^{\infty}a(k)k^{- s},\quad\quad Res > 1,\quad k^{- s}: = \exp\left( - s\log k \right).
\tag{1.2}\label{eq:1.2}
\end{equation}

The logarithm is real since \(k > 0.\) Absolute convergence follows by
separating the first \(K\) terms and bounding the remainder by
\(A\sum_{k > K}^{}k^{- Res}.\) All samples below therefore lie in the
absolute-convergence region.

The observations are a vector
\(y = \left( y_{0},\ldots,y_{K - 1} \right) \in \mathbb{C}^{K}\) that
satisfying

\begin{equation}
y_{j} = F_{a}(2q + 2j) + e_{j},\quad\quad\left| e_{j} \right| \leq \delta,\quad\quad 0 \leq j < K,
\tag{1.3}\label{eq:1.3}
\end{equation}

where \(\delta \geq 0\) is a known absolute error bound and no
probability distribution is assumed for the errors.

The classical theory of Dirichlet series provides the setting for this
recovery problem {[}1{]}. The observation block is also a finite moment
sequence, connecting the problem with Hausdorff's work {[}2{]}. Its
finite-dimensional inversion uses Lagrange interpolation {[}3{]} and
inverse Vandermonde matrices {[}4{]}. The question of the smallest
attainable worst-case error belongs to deterministic optimal recovery
{[}5{]}.

Through the studies, Choulli {[}6{]} worked on Dirichlet-coefficient
recovery in connection with inverse heat problems. Gerth et.al {[}7{]}
also treat stability and truncated observations for the Hausdorff moment
problem in Hilbert spaces. Foucart and Liao {[}8{]} study recovery with
Euclidean bounds on observational errors. Here the errors are bounded
separately in each coordinate, the tail coefficients are bounded
individually and the initial block is unrestricted. These assumptions
allow the complete feasible set of the target coefficient to be
determined for every observed vector.

Here, we define the degree-\((K - 1)\) polynomial

\begin{equation}
p_{n,K,q}(t) = n^{2q}\prod_{\substack{1 \leq m \leq K \\ m \neq n}}^{}\frac{t - m^{- 2}}{n^{- 2} - m^{- 2}} = \sum_{j = 0}^{K - 1}c_{j}t^{j}.
\tag{1.4}\label{eq:1.4}
\end{equation}

We clear products equal \(1,\) including the case \(K = n = 1.\) Its
coefficients \(c_{j}\) are real rational numbers. The estimator and its
tail weights are

\begin{equation}
{\widehat{a}}_{n}(y) = \sum_{j = 0}^{K - 1}c_{j}y_{j},\quad\quad w_{k} = k^{- 2q}p_{n,K,q}\left( k^{- 2} \right)\quad(k > K).
\tag{1.5}\label{eq:1.5}
\end{equation}

The polynomial is a scaled Lagrange cardinal polynomial which equals
\(n^{2q}\) at \(n^{- 2}\) and vanishes at the other first \(K\) nodes.

We set the tail constant, the noise amplification and their combined
radius by

\begin{equation}
B_{n,K,q} = \sum_{k > K}^{}\left| w_{k} \right|,\quad\quad\Lambda_{n,K,q} = \sum_{j = 0}^{K - 1}\left| c_{j} \right|,\quad\quad E_{n,K,q} = AB_{n,K,q} + \delta\Lambda_{n,K,q}.
\tag{1.6}\label{eq:1.6}
\end{equation}

The tail sum converges because \(p_{n,K,q}\) is bounded on
\(\lbrack 0,1\rbrack\) and \(2q > 1.\) Our first result gives the full
uncertainty set and its optimal center.

\textbf{Theorem 1.1 (Exact uncertainty disk and minimax recovery).} For
every \(y \in \mathbb{C}^{K},\) the set of coefficient values compatible
with (1.1) and (1.3) is exactly

\begin{equation}
\mathcal{F}_{n}(y): = \{ a(n):(a,e)\text{ is compatible with }y\} = \{ z\in \mathbb{C}:\left| z - {\widehat{a}}_{n}(y) \right| \leq E_{n,K,q}\}.
\tag{1.7}\label{eq:1.7}
\end{equation}

Consequently, among all maps
\(\Phi:\mathbb{C}^{K}\rightarrow \mathbb{C},\) including nonlinear maps,

\begin{equation}
\inf_{\Phi}\ \sup_{\substack{a \in \mathcal{A}_{K}(A) \\ \left| e_{j} \right| \leq \delta}}\left| \Phi(\left( F_{a}(2q + 2j) + e_{j} \right)_{j = 0}^{K - 1}) - a(n) \right| = E_{n,K,q}.
\tag{1.8}\label{eq:1.8}
\end{equation}

The estimator (1.5) attains this value and is the unique optimal center
at each fixed \(y.\)

To state the second result, we put

\begin{equation}
C_{n,K,q} = \frac{2n^{2q + 2K}}{(K - n)!(K + n)!},
\tag{1.9}\label{eq:1.9}
\end{equation}

where \(0! = 1.\) The notation \(u_{K} \sim v_{K}\) means
\(u_{K}/v_{K} \rightarrow 1\) as \(K \rightarrow \infty.\)

\textbf{Theorem 1.2 (Sharp asymptotic tail scale).} For fixed integers
\(n,q \geq 1,\)

\begin{equation}
B_{n,K,q} \sim C_{n,K,q}\,\frac{3^{q - 1/2}}{2}\Gamma\left( q - \frac{1}{2} \right)K^{\frac{3}{2}(1 - 2q)}\quad\quad(K \rightarrow \infty)
\tag{1.10}\label{eq:1.10}
\end{equation}

where \(\Gamma(v) = \int_{0}^{\infty}u^{v - 1}e^{- u}\, du\) for
\(v > 0.\) The dominant tail indices are of order \(K^{3/2},\) in the
normalized-mass sense specified in the proof.

The interpolation and norm principles used here are classical and the
results give their explicit consequences for reciprocal-square nodes and
the stated uncertainty class.

\hypertarget{interpolation-identities-and-explicit-constants}{%
\section{Interpolation identities and explicit
constants}\label{interpolation-identities-and-explicit-constants}}

The samples in (1.3), before noise is added, are moments of the finite
complex measure \(\sum_{k \geq 1}^{}a(k)x_{k}^{q}\delta_{x_{k}}\) on
\(\lbrack 0,1\rbrack,\) where \(\delta_{x_{k}}\) is unit point mass at
\(x_{k}.\) The measure has finite total variation since
\(\sum_{k \geq 1}^{}\left| a(k) \right|k^{- 2q} < \infty.\) The cardinal
values of (1.4) give the following identity.

\textbf{Lemma 2.1 (Exact error identity).} For every admissible pair
\((a,e)\) in (1.1)--(1.3),

\begin{equation}
{\widehat{a}}_{n}(y) - a(n) = \sum_{k > K}^{}w_{k}a(k) + \sum_{j = 0}^{K - 1}c_{j}e_{j}
\tag{2.1}\label{eq:2.1}
\end{equation}

The first sum is absolutely convergent. The estimator is the unique
complex-linear rule that is exact for every sequence supported on
\(\{ 1,\ldots,K\}\) with zero observational error.

\emph{Proof.} we insert (1.3) into (1.5). There are finitely many \(j\)
and each Dirichlet series is absolutely convergent, so we have

\begin{equation}
\sum_{j = 0}^{K - 1}c_{j}F_{a}(2q + 2j) = \sum_{k = 1}^{\infty}a(k)k^{- 2q}p_{n,K,q}\left( k^{- 2} \right)
\tag{2.2}\label{eq:2.2}
\end{equation}

For \(k \leq K,\) the factor multiplying \(a(k)\) is \(1\) if \(k = n\)
and \(0\) otherwise. This proves (2.1). The polynomial is bounded on
\(\lbrack 0,1\rbrack\) and the coefficient tail is bounded by \(A,\)
therefore the tail sum is absolutely convergent. Any exact
complex-linear rule has coefficients \(d_{j}\) whose polynomial
\(d(t) = \sum_{j < K}^{}d_{j}t^{j}\) satisfies
\(k^{- 2q}d\left( k^{- 2} \right) = \mathbf{1}_{\{ k = n\}}\) for
\(k \leq K,\) as follows by testing sequences with a single nonzero
entry. Thus \(d - p_{n,K,q}\) has \(K\) distinct zeros and degree at
most \(K - 1\) and must vanish identically. The symbol
\(\mathbf{1}_{\{ k = n\}}\) is \(1\) when \(k = n\) and \(0\) otherwise.

We next compute the constants in (1.6). In addition to (1.9), we define
the elementary tail bound

\begin{equation}
T_{n,K,q} = \frac{C_{n,K,q}K^{1 - 2q}}{2q - 1}
\tag{2.3}\label{eq:2.3}
\end{equation}

\textbf{Proposition 2.2 (Closed constants and certified tail
evaluation).} For \(q \geq 1\) and \(K \geq n \geq 1,\)

\begin{equation}
\left| p_{n,K,q}(0) \right| = C_{n,K,q},
\tag{2.4}\label{eq:2.4}
\end{equation}

\begin{equation}
\Lambda_{n,K,q} = C_{n,K,q}\prod_{\substack{1 \leq m \leq K \\ m \neq n}}^{}\left( 1 + m^{2} \right),
\tag{2.5}\label{eq:2.5}
\end{equation}

\begin{equation}
B_{n,K,q} = C_{n,K,q}\sum_{k > K}^{}k^{- 2q}\prod_{\substack{1 \leq m \leq K \\ m \neq n}}^{}\left( 1 - \frac{m^{2}}{k^{2}} \right) \leq T_{n,K,q}.
\tag{2.6}\label{eq:2.6}
\end{equation}

For an integer \(J > K,\) let us \(B_{n,K,q}^{\lbrack J\rbrack}\) denote
the sum in (2.6) truncated at \(k = J,\) including its prefactor. Then

\begin{equation}
0 \leq B_{n,K,q} - B_{n,K,q}^{\lbrack J\rbrack} \leq \frac{C_{n,K,q}J^{1 - 2q}}{2q - 1}.
\tag{2.7}\label{eq:2.7}
\end{equation}

Consequently finite rational calculations enclose the tail constant and,
for rational \(A\) and \(\delta,\) the exact minimax radius.

\emph{Proof.} At zero, each factor in (1.4) has absolute value
\(n^{2}/\left| m^{2} - n^{2} \right|.\) Its denominator product
satisfies

\begin{equation}
\prod_{\substack{m \leq K \\ m \neq n}}^{}\left| m^{2} - n^{2} \right| = (n - 1)!(K - n)!\,\frac{(K + n)!}{n!\, 2n} = \frac{(K - n)!(K + n)!}{2n^{2}}
\tag{2.8}\label{eq:2.8}
\end{equation}

Multiplying by \(n^{2q + 2K - 2}\) proves (2.4). Exactly the factors
with \(m > n\) are negative at zero, so its sign is \(( - 1)^{K - n}.\)

Factoring out the constant term gives us

\begin{equation}
p_{n,K,q}(t) = p_{n,K,q}(0)\prod_{\substack{m \leq K \\ m \neq n}}^{}\left( 1 - m^{2}t \right).
\tag{2.9}\label{eq:2.9}
\end{equation}

The product coefficients alternate in sign which their absolute values
are the nonnegative elementary symmetric sums of the numbers \(m^{2}.\)
Their sum equals the product evaluated at \(t = - 1,\) proving (2.5) and
this specializes the coefficient-norm mechanism in Gautschi's
inverse-Vandermonde analysis {[}4{]}.

For \(k > K,\) each factor \(1 - m^{2}/k^{2}\) lies between zero and
one. Thus all \(w_{k}\) have sign \(( - 1)^{K - n}\) and their absolute
values give (2.6). Finally, integral comparison yields

\begin{equation}
\sum_{k > J}^{}k^{- 2q} \leq \int_{J}^{\infty}t^{- 2q}\, dt = \frac{J^{1 - 2q}}{2q - 1}.
\tag{2.10}\label{eq:2.10}
\end{equation}

Taking \(J = K\) proves the bound by \(T_{n,K,q};\) taking arbitrary
\(J > K\) proves Eq (2.7).

For instance, \(n = K = 2\) and \(q = 1\) give
\(p_{2,2,1}(t) = \frac{16}{3}(1 - t),\)
\({\widehat{a}}_{2} = \frac{16}{3}\left( y_{0} - y_{1} \right)\) and
\(\Lambda_{2,2,1} = 32/3.\) The exact radius is

\begin{equation}
E_{2,2,1} = \frac{16A}{3}\sum_{k = 3}^{\infty}\left( k^{- 2} - k^{- 4} \right) + \frac{32\delta}{3}.
\tag{2.11}\label{eq:2.11}
\end{equation}

The positivity of this tail sum shows us why cancellation estimates for
a special coefficient sequence would not suffice to establish the
worst-case result.

\hypertarget{proof-of-theorem-1.1}{%
\section{Proof of Theorem 1.1}\label{proof-of-theorem-1.1}}

\emph{Proof.} Lemma 2.1 and the triangle inequality give inclusion in
the disk. For proving the reverse inclusion, we preserve the entire
observed vector while varying the tail and noise.

We fix an arbitrary admissible tail \(\left( a(k) \right)_{k > K}\) and
errors \(\left( e_{j} \right)_{j < K}.\) The first \(K\) coefficients
msut satisfy

\begin{equation}
\sum_{k = 1}^{K}a(k)k^{- 2q}x_{k}^{j} = y_{j} - e_{j} - \sum_{k > K}^{}a(k)k^{- 2q}x_{k}^{j},\quad\quad 0 \leq j < K.
\tag{3.1}\label{eq:3.1}
\end{equation}

This square system is invertible. A vector \(h\) in its kernel satisfies
\(\sum_{k \leq K}^{}h(k)k^{- 2q}x_{k}^{j} = 0\) for all \(j < K.\)
Taking the linear combination specified by \(p_{m,K,q}\) gives
\(h(m) = 0\) for each \(m \leq K.\) The kernel is therefore zero. The
weights \(w_{k}\) are nonzero for \(k > K,\) since their arguments are
distinct from the polynomial roots. For any \(v\in \mathbb{C}\) with
\(|v| \leq 1,\) we choose

\begin{equation}
a(k) = Av\frac{\overline{w_{k}}}{\left| w_{k} \right|}\ (k > K),\quad\quad e_{j} = \delta v\frac{\overline{c_{j}}}{\left| c_{j} \right|}\ \left( c_{j} \neq 0 \right)
\tag{3.2}\label{eq:3.2}
\end{equation}

and set \(e_{j} = 0\) if \(c_{j} = 0.\) The bar denotes complex
conjugation. We solve (3.1) for the first block. Equation (2.1) then
gives \(a(n) = {\widehat{a}}_{n}(y) - vE_{n,K,q},\) covering the entire
disk.

For any proposed value \(z_{0},\) the largest distance from \(z_{0}\) to
this disk is
\(\left| z_{0} - {\widehat{a}}_{n}(y) \right| + E_{n,K,q}.\) The upper
bound is the triangle inequality (equality is attained at the boundary
point opposite \(z_{0},\) with any boundary point used if \(z_{0}\) is
the center). The unique minimizing choice is
\(z_{0} = {\widehat{a}}_{n}(y).\) The radius is independent of \(y\)
which proving the global statement (1.8).

For real coefficients, observations and errors, the same construction
yields the interval
\(\left\lbrack {\widehat{a}}_{n}(y) - E_{n,K,q},{\widehat{a}}_{n}(y) + E_{n,K,q} \right\rbrack.\)
The first-block freedom is crucial which imposing
\(\left| a(k) \right| \leq A\) also for \(k \leq K\) can shrink the
feasible set. The upper bound remains valid on that smaller class, but
the minimax equality is not asserted there.

\hypertarget{proof-of-theorem-1.2}{%
\section{Proof of Theorem 1.2}\label{proof-of-theorem-1.2}}

\emph{Proof.} Here we set \(h_{K} = K^{3/2}\) and
\(P_{K}(k) = \prod_{1 \leq m \leq K,\, m \neq n}^{}\left( 1 - m^{2}/k^{2} \right)\)
for \(k > K.\) By defining a nonnegative step function \(f_{K}\) to be
zero on \(\left( 0,K/h_{K} \right)\),we can satisfy

\begin{equation}
f_{K}(t) = \left( k/h_{K} \right)^{- 2q}P_{K}(k)\quad\text{for }(k - 1)/h_{K} \leq t < k/h_{K},\quad k > K.
\tag{4.1}\label{eq:4.1}
\end{equation}

Endpoint values do not affect its integral. Direct integration gives us

\begin{equation}
\int_{0}^{\infty}f_{K}(t)\, dt = h_{K}^{2q - 1}\frac{B_{n,K,q}}{C_{n,K,q}}.
\tag{4.2}\label{eq:4.2}
\end{equation}

We fix \(t > 0\ and\ \)for large \(K,\) its interval index \(k\) exceeds
\(K\) and \(k/h_{K} \rightarrow t.\) By Putting
\(S_{K} = \sum_{1 \leq m \leq K,\, m \neq n}^{}m^{2},\) the
sum-of-squares formula gives \(S_{K}/K^{3} \rightarrow 1/3,\) while
\(\max_{m \leq K}m^{2}/k^{2} \rightarrow 0.\) For \(0 \leq u \leq 1/2,\)
\(\left| \log(1 - u) + u \right| \leq u^{2},\) because
\(\sum_{r \geq 2}^{}u^{r}/r \leq u^{2}/\left\lbrack 2(1 - u) \right\rbrack.\)
Hence

\begin{equation}
\log P_{K}(k) = - \frac{S_{K}}{k^{2}} + O\left( \frac{\sum_{m = 1}^{K}m^{4}}{k^{4}} \right) \rightarrow - \frac{1}{3t^{2}}.
\tag{4.3}\label{eq:4.3}
\end{equation}

Here the error is \(O\left( K^{- 1} \right)\) for fixed \(t,\) since
\(\sum m^{4} \leq K^{5}\) and \(k\) is of order \(K^{3/2}.\) The \(O\)
notation means an absolute bound by a constant multiple of its argument.
In this step, the constant may depend on fixed \(t.\) Thus
\(f_{K}(t) \rightarrow t^{- 2q}\exp\left\lbrack - 1/\left( 3t^{2} \right) \right\rbrack.\)
For domination, we choose \(K\) large enough that
\(S_{K} \geq K^{3}/4.\) The inequality \(\log(1 - u) \leq - u\) gives
\(P_{K}(k) \leq exp\left\lbrack - S_{K}/k^{2} \right\rbrack.\) On a
nonzero interval of (4.1), \(t \leq k/h_{K} \leq 2t\) because
\(k \geq 2.\) Accordingly,

\begin{equation}
0 \leq f_{K}(t) \leq t^{- 2q}\exp\left( - \frac{1}{16t^{2}} \right).
\tag{4.4}\label{eq:4.4}
\end{equation}

This majorant is integrable at zero by exponential decay and at infinity
because \(2q > 1.\) Dominated convergence therefore applies to (4.2).
Substituting \(u = 1/\left( 3t^{2} \right)\) in the limit integral
yields

\begin{equation}
\int_{0}^{\infty}t^{- 2q}e^{- 1/\left( 3t^{2} \right)}\, dt = \frac{3^{q - 1/2}}{2}\int_{0}^{\infty}u^{q - 3/2}e^{- u}\, du,
\tag{4.5}\label{eq:4.5}
\end{equation}

which proves (1.10). Restricting the normalized sum to
\(\alpha h_{K} \leq k \leq \beta h_{K},\) for fixed
\(0 < \alpha < \beta < \infty,\) gives the limit integral over
\(\lbrack\alpha,\beta\rbrack\) by the same dominated-convergence
argument. Endpoint intervals have vanishing width and bounded height on
\(\lbrack\alpha/2,2\beta\rbrack.\) Dividing by the positive total limit
in (4.5) gives the limiting fraction of the tail mass in this interval.
The remaining limiting mass tends to zero as \(\alpha \downarrow 0\) and
\(\beta \uparrow \infty,\) establishing the stated dominant scale. This
asymptotic holds for fixed \(n\) and \(q,\) not uniformly for target
indices growing with \(K.\) At finite \(K,\) (2.7) remains a rigorous
enclosure.

\hypertarget{selecting-the-sample-count}{%
\section{Selecting the sample
count}\label{selecting-the-sample-count}}

Proposition 2.2 yields the computable certificate

\begin{equation}
\left| {\widehat{a}}_{n}(y) - a(n) \right| \leq E_{n,K,q} \leq U_{n,K,q}: = AT_{n,K,q} + \delta\Lambda_{n,K,q}.
\tag{5.1}\label{eq:5.1}
\end{equation}

The first term decreases rapidly for large \(K,\) but the second
generally increases. Their behavior determines how many samples should
be used in this estimator.

\textbf{Proposition 5.1 (Noise amplification).} For fixed \(n,q\) and
\(K \geq n,\)

\begin{equation}
\frac{C_{n,K + 1,q}}{C_{n,K,q}} = \frac{n^{2}}{(K + 1)^{2} - n^{2}},\quad\quad\frac{\Lambda_{n,K + 1,q}}{\Lambda_{n,K,q}} = \frac{n^{2}\{ 1 + (K + 1)^{2}\}}{(K + 1)^{2} - n^{2}}.
\tag{5.2}\label{eq:5.2}
\end{equation}

The sequence \(\Lambda_{n,K,q}\) is strictly increasing. With
\(P_{\infty}: = \prod_{m = 1}^{\infty}\left( 1 + m^{- 2} \right),\)

\begin{equation}
\lim_{K \rightarrow \infty}\frac{\Lambda_{n,K,q}}{n^{2K}} = \frac{2n^{2q}}{1 + n^{2}}P_{\infty}.
\tag{5.3}\label{eq:5.3}
\end{equation}

It is bounded for \(n = 1\) and grows exponentially for each fixed
\(n \geq 2.\)

\emph{Proof.} Dividing successive instances of Eqs (1.9) and (2.5) gives
(5.2). The second ratio exceeds one which its denominator is positive
and numerator minus denominator is
\(\left( n^{2} - 1 \right)(K + 1)^{2} + 2n^{2} > 0.\)

Factoring \(m^{2}\) from each factor of (2.5) gives

\begin{equation}
\frac{\Lambda_{n,K,q}}{n^{2K}} = \frac{2n^{2q}}{1 + n^{2}}\frac{(K!)^{2}}{(K - n)!(K + n)!}\prod_{m = 1}^{K}\left( 1 + m^{- 2} \right).
\tag{5.4}\label{eq:5.4}
\end{equation}

The factorial ratio is \(\prod_{r = 0}^{n - 1}(K - r)/(K + r + 1),\)
which is at most one and tends to one. Since
\(0 < log\left( 1 + m^{- 2} \right) \leq m^{- 2}\) and
\(\sum_{m \geq 1}^{}m^{- 2} \leq 2,\) the product \(P_{\infty}\) exists,
is positive and is at most \(e^{2}.\) Taking limits proves (5.3).

\textbf{Theorem 5.2 (A terminating global certificate).} We fix
\(n \geq 2,\) \(q \geq 1,\) \(A > 0\), \(\delta > 0\ and\ \)examine
\(K = n,n + 1,\ldots,\) retaining the smallest certificate \(U_{best}\)
already obtained from (5.1). We stop at the first unexamined integer
\(K_{0}\) for which

\begin{equation}
\delta\Lambda_{n,K_{0},q} \geq U_{best}.
\tag{5.5}\label{eq:5.5}
\end{equation}

The procedure terminates and a retained minimizer minimizes
\(U_{n,K,q}\) over all integers \(K \geq n.\) Ties may be resolved by
retaining the smallest examined \(K.\) The result concerns \(U,\) not
necessarily the exact radius \(E.\)

\emph{Proof.} The value at \(K = n\) is finite, whereas Proposition 5.1
gives us \(\delta\Lambda_{n,K,q} \rightarrow \infty.\) The stopping
condition must eventually hold. For every \(K \geq K_{0},\) monotonicity
gives
\(U_{n,K,q} \geq \delta\Lambda_{n,K,q} \geq \delta\Lambda_{n,K_{0},q} \geq U_{best}.\)
All smaller integers were explicitly compared, proving global minimality
of the retained value. When \(K\) varies, the tail class (1.1) also
varies. This theorem minimizes a family of valid certificates, not
minimax errors on one unchanged class. If \(\left| a(k) \right| \leq A\)
for all \(k,\) every certificate applies to that same sequence and the
minimax equality still refers to the larger class at each fixed \(K.\)

\textbf{Corollary 5.3 (Decreasing-noise limit and integer recovery).} We
fix \(n,q,A\ and\ \)for \(0 < \delta < e^{- 1}\) set
\(K(\delta) = max\{ n,\lfloor\sqrt{\log(1/\delta)}\rfloor\}.\) Then
\(U_{n,K(\delta),q} \rightarrow 0\) as \(\delta \downarrow 0.\) For a
real integer target coefficient, a certificate \(U_{n,K,q} < 1/2\)
guarantees exact recovery by rounding \(Re{\widehat{a}}_{n}(y)\) to the
nearest integer.

\emph{Proof.} The first ratio in (5.2) tends to zero, so
\(C_{n,K,q} \rightarrow 0\) and consequently
\(T_{n,K,q} \rightarrow 0.\) Equation (5.4) bounds the noise term by
\(2e^{2}n^{2q}\delta n^{2K}/\left( 1 + n^{2} \right).\) With
\(L = log(1/\delta),\) its varying part is at most
\(\exp\left\lbrack - L + 2\left( \sqrt{L} + n \right)\log n \right\rbrack,\)
which tends to zero. This proves the limit. Finally, the absolute
real-part error is bounded by the modulus of the complex error and
distinct integers are separated by at least one.

\section*{Dedication}
To the memory of Bernhard Riemann, whose mathematics has profoundly
shaped the way we think. I dedicate this work to him with deep gratitude
and the hope of contributing something worthy of the questions he left
us.

\hypertarget{references}{%
\section*{References}\label{references}}

{[}1{]} ``The General Theory of Dirichlet's Series \textbar{} Cambridge
University Press \& Assessment.'' Available:
\url{https://www.cambridge.org/it/universitypress/subjects/mathematics/historical-mathematical-texts/general-theory-dirichlets-series}

{[}2{]} F. Hausdorff, ``Summationsmethoden und Momentfolgen. I,''
\emph{Mathematische Zeitschrift}, vol. 9, no. 1--2, pp. 74--109, Mar.
1921, doi: 10.1007/BF01378337/METRICS.

{[}3{]} J. P. Berrut and L. N. Trefethen, ``Barycentric Lagrange
Interpolation,'' \url{https://doi.org/10.1137/S0036144502417715}, vol.
46, no. 3, pp. 501--517, Aug. 2006, doi: 10.1137/S0036144502417715.

{[}4{]} W. Gautschi, ``On inverses of Vandermonde and confluent
Vandermonde matrices,'' \emph{Numer. Math. (Heidelb).}, vol. 4, no. 1,
pp. 117--123, Dec. 1962, doi: 10.1007/BF01386302/METRICS.

{[}5{]} C. A. Micchelli and T. J. Rivlin, ``A Survey of Optimal
Recovery,'' \emph{Optimal Estimation in Approximation Theory}, pp.
1--54, 1977, doi: 10.1007/978-1-4684-2388-4\_1.

{[}6{]} M. Choulli, ``Various stability estimates for the problem of
determining an initial heat distribution from a single measurement,''
Dec. 2015. Available: https://arxiv.org/pdf/1512.07421

{[}7{]} D. Gerth, B. Hofmann, C. Hofmann, and S. Kindermann, ``The
Hausdorff Moment Problem in the light of ill-posedness of type I,''
\emph{Eurasian Journal of Mathematical and Computer Applications}, vol.
9, no. 2, pp. 57--87, 2021, doi: 10.32523/2306-6172-2021-9-2-12-38.

{[}8{]} S. Foucart and C. Liao, ``Optimal Recovery from Inaccurate Data
in Hilbert Spaces: Regularize, But What of the Parameter?,''
\emph{Constructive Approximation 2022 57:2}, vol. 57, no. 2, pp.
489--520, Oct. 2022, doi: 10.1007/S00365-022-09590-5.

~

Polytechnic of Turin

\href{mailto:Payam.danesh@polito.it}{\nolinkurl{Payam.danesh@polito.it}}

\end{document}